\documentstyle{article}
\def\Bbb{\bf}
\newcommand{\IM}{\mathop{\rm Im}}
\newcommand{\Gr}{\mbox{\rm Gr}}

\newcommand{\rank}{\mathop{\rm rank}}

\newcommand{\tr}{\mathop{\rm tr}}
\newcommand{\Char}{ {\mbox{\rm  \, char}^*\, } }
\newcommand{\Span}{ {\mbox{\rm Span}} }
\newcommand{\Stab}{ {\mbox{\rm Stab}} }
\newcommand{\diag}{\mathop{\rm diag}}
\newtheorem{thm}{Theorem}
\newtheorem{lem}[thm]{Lemma}
\newtheorem{cor}[thm]{Corollary}
\newtheorem{prop}[thm]{Proposition}
\newtheorem{conj}[thm]{Conjecture}
\newcommand{\eqdef}{\stackrel{\rm def}{=}}
\newcommand{\CC}{{\Bbb C}}

\newcommand{\ZZ}{{\Bbb Z}}
\newcommand{\NN}{{\Bbb N}}
\newcommand{\PP}{{\Bbb P}}
\newcommand{\al}{{\alpha}}
\newcommand{\be}{{\beta}}
\newcommand{\Del}{{\Delta}}
\newcommand{\tDel}{\widetilde{\Del}}
\newcommand{\lam}{ \lambda }
\newcommand{\Lam}{ \Lambda }
\newcommand{\om}{ \varpi }
\newcommand{\FF}{ { \cal F } }
\newcommand{\FFB}{ \FF^B }
\newcommand{\MB}{ M^B }
\newcommand{\Mt}{ M^{\mbox{\tiny tensor}} }
\newcommand{\FFD}{ \FF_D }

\newcommand{\FFBD}{ \FF^B_D }
\newcommand{\LL}{ {\cal L} }

\newcommand{\OO}{ {\cal O} }
\newcommand{\SS}{ {\cal S} }
\newcommand{\II}{ {\cal I} }

\newcommand{\ii}{ {\bf i} }
\newcommand{\jj}{ {\bf j} }
\newcommand{\mm}{ {\bf m} }
\newcommand{\ww}{ {\bf w} }
\newcommand{\LLmm}{ \LL_{\mm} }
\newcommand{\Zii}{ Z_{\ii}  }
\newcommand{\zii}{ z_{\ii}  }

\newcommand{\Dwj}{ D_{\ww\jj}  }
\newcommand{\ziiq}{ z_{\ii}^{\tiny quo}  }
\newcommand{\Fii}{ \Zii^{\tiny fib}  }
\newcommand{\Qii}{ \Zii^{\tiny quo}  }
\newcommand{\Oii}{ \Zii^{\tiny orb}  }
\newcommand{\Gii}{ \Zii^{\tiny emb}  }
\newcommand{\Dii}{ D_{\ii}  }
\newcommand{\Dtii}{ D^+_{\ii}  }
\newcommand{\Ph}{ \widehat{P} }

\newcommand{\tw}{ \tilde{w} }
\newcommand{\tu}{ \tilde{u} }
\newcommand{\tww}{ \widetilde{\ww} }
\newcommand{\tmm}{ \widetilde{\mm} }
\newcommand{\tD}{ \widetilde{D} }
\newcommand{\mtimes}{ \mathop{\times} }
\newcommand{\lelt}{ \stackrel{e\!l\!t}{<} }
\newcommand{\leqcomp}{ \stackrel{c\!o\!m\!p}{\leq} }
\newcommand{\leqlex}{ \stackrel{l\!e\!x}{\leq} }

\title{Bott-Samelson Varieties \\
 and Configuration Spaces}
\author{Peter Magyar}
\date{October, 1996}
\begin{document}

\maketitle

\begin{center} {\bf Abstract} \end{center}

{\small
\noindent
The Bott-Samelson varieties $Z$ are a powerful tool
in the representation theory and geometry of
a reductive group $G$.  We give a new construction of
$Z$ as the closure of a $B$-orbit in 
a product of flag varieties $(G/B)^l$.
This also gives an embedding of the projective
coordinate ring of the variety into the function 
ring of a Borel subgroup: $\CC[Z] \subset \CC[B]$.

In the case of the general linear group $G = GL(n)$, 
this identifies $Z$ as a configuration variety of
multiple flags subject to certain inclusion conditions,
controlled by the combinatorics
of braid diagrams and generalized Young diagrams. 
The natural mapping $Z \rightarrow G/B$ compactifies
the matrix factorizations of Berenstein, 
Fomin and Zelevinsky \cite{BFZ}.
As an application, we give a geometric proof of the
theorem of Kraskiewicz and Pragacz \cite{KP} that
Schubert polynomials are characters of Schubert modules.

Our work leads on the one hand
to a Demazure character formula for
Schubert polynomials and 
other generalized Schur functions,
and on the other hand to a Standard Monomial
Theory for Bott-Samelson varieties.
All our results remain valid in arbitrary
characteristic and over $\ZZ$.
}

\vspace{1em}
\noindent
{\large \bf Introduction}
\\[1em]
The Bott-Samelson varieties
are an important geometric
tool in the theory of a reductive 
algebraic group 
(or complex Lie group) $G$.
Defined in \cite{BS}, 
they were exploited by Demazure
\cite{Dem1} to analyze the flag variety $G/B$,
its singular cohomology ring
$H^{\cdot}(G/B,\CC)$ (the Schubert calculus),
and its projective coordinate ring 
$\CC[G/B]$.  Since the irreducible representations 
of $G$ are embedded in the coordinate ring,
Demazure was able to obtain an
iterative character formula \cite{Dem2}
for these representations.

Bott-Samelson
varieties are so useful because they 
``factor'' the flag variety
into a ``product'' of projective lines.
More precisely, they are iterated 
$\PP^1$-fibrations
and each has a natural, birational map to 
$G/B$.
The Schubert subvarieties themselves
lift to iterated $\PP^1$-fibrations
under this map.
The combinatorics of Weyl groups enters
the picture because a given $G/B$
can be ``factored'' in many 
ways, indexed by sequences 
$\ii = (i_1, i_2, \ldots, i_N)$
such that $w_0 = s_{i_1} s_{i_2} \cdots s_{i_N}$
is a reduced decompostion of the longest 
Weyl group element $w_0$ into simple reflections.

The Bott-Samelson variety $\Zii$
is usually defined as a quotient:
$$
\Zii \eqdef 
(P_{i_1} \mtimes \cdots \mtimes P_{i_N}) / B^N,
$$
where $P_i$ are minimal parabolic
subgroups, $B \subset P_i \subset G$, 
and $B^N$ acts freely on the right
of $P_{i_1} \mtimes \cdots \mtimes P_{i_N}$ by 
$$
(p_1,\ldots,p_N) \cdot (b_1,\ldots,b_N) 
= (p_1 b_1, b_1^{-1} p_2 b_2, 
\ldots, b_{N-1}^{-1} p_N b_N).
$$
The natural map to the flag variety 
is given by multiplication:
$(p_1,\ldots,p_N) \mapsto$ \\
$p_1 p_2 \cdots p_N B\in G/B$.

In this paper, we first
give a dual construction of
$\Zii$ as a subvariety rather than a quotient.
It is the closure of a $B$-orbit 
inside a product of flag varieties:
$$
\Zii \cong \overline{B\cdot(s_{i_1} B, 
s_{i_1}\! s_{i_2} B, \ldots, w_0 B)} \subset (G/B)^N,
$$
where $B$ acts diagonally on $(G/B)^N$.
Our constructions are partly inspired by
Fulton's work \cite{F}, Ch. 10.3.

In the case $G=GL(n)$ or $SL(n)$, 
this translates into
an expression for $\Zii$ as a 
``multiple Schubert variety'':
configurations of many linear spaces in 
$\CC^n$ subject to certain inclusions
involving a test flag.
For example, for $G =GL(3)$,
$\ii = 212$, and the test flag 
$\CC^1 \subset\CC^2 \subset \CC^3$,
we get
$$
\Zii = \{(V_1, V_2, V_2') \in 
\Gr(1,\CC^3) \times \Gr(2,\CC^3)^2 \mid
V_2 \supset V_1 \subset V_2' \supset \CC^1 \}.
$$
The natural birational map onto the flag variety is
given by the projection
$(V_1,V_2,V_2') \mapsto (V_1,V_2)$.
For $GL(n)$, the combinatorics of 
such configuration varieties
is controlled by certain
generalized Young diagrams 
\cite{MNW}, \cite{MFour}, \cite{RS1}, \cite{RSNew}; 
or equivalently by the wiring diagrams and 
chamber sets of Berenstein, Fomin, and Zelevinsky
\cite{BFZ}, \cite{LZ}.

Secondly, we study more general
configuration varieties,
which are also closures of $B$-orbits in products
of $G/B$.  These varieties are governed by
similar combinatorics, are desingularized by
the Bott-Samelson varieties, and  
include the flag and Schubert varieties.

Thirdly, we turn to the Borel-Weil 
theory of Bott-Samelson varieties.  
Our embedding of $\Zii$ leads to
an embedding of its projective coordinate ring
into the regular functions on a Borel
subgroup:
$$
\CC[\Zii] \subset \CC[B].
$$
That is, the space of
sections of effective line bundles
on $\Zii$ can be realized in terms of certain
polynomials on $B$.
(Here we use a vanishing theorem of W. van der Kallen
\cite{MNW}.)

For $G = GL(n)$, the space of sections
becomes a certain
generalized Schur module 
(\cite{ABW}, \cite{Ro}, \cite{RS1}, \cite{RS3}, \cite{RSNew})
spanned by products
of minors in the polynomial ring $\CC[x_{ij}]_{i<j}$.
Here, the bitableaux of Desarmenien, Kung, and Rota
\cite{DKR} 
(c.f. \cite{LLT}), 
give the appropriate combinatorial formalism.
A result of our construction 
is a Demazure character
formula for these generalized Schur modules.
Conversely, we get a standard monomial
basis for the space of sections, which we pursue
in our paper \cite{LM}.

Fourthly, we apply our results 
to the Schubert modules of Kraskiewicz 
and Pragacz \cite{KP}. The characters 
of these modules are the Schubert polynomials,
special algebraic representatives 
of the Schubert classes
in the singular cohomology ring of $G/B$.
Why the Schubert polynomials should
appear as characters of $B$-modules 
remains a mystery, but our theory does lead 
(as suggested by a manuscript of V. Reiner
and M. Shimozono) to
a new proof of Kraskiewicz and Pragacz's theorem.
Our Demazure formula applies to these polynomials,
and is basically different from 
the usual recurrence defining 
them.  The combinatorics of this 
formula are examined
in our paper \cite{MFour}.

To avoid intimidating terminology, we work over the
base field $\CC$ of complex numbers.  The alert reader will
note, however, that all our arguments remain
valid without change over an algebraically closed
field of arbitrary characteristic and over the
integers.

{\bf Note.} The geometry of a general reductive $G$ is 
largely confined to Sec 1 and 2.
Those interested mainly in the combinatorial 
applications associated to $G=GL(n)$ 
may begin reading at Sec 3.
\\[1em]
{\bf Acknowledgements.}  The author would like to thank
Victor Reiner, Mark Shimozono, and Bill Fulton
for numerous helpful suggestions 
and for making available their unpublished work.
\\[1em]
{\small
{\bf Contents.} {\bf 1.} Bott-Samelson varieties \ 
{\bf 1.1}  Three constructions \
{\bf 1.2}  Isomorphism theorem \
{\bf 1.3}  Open cells \ \
{\bf 2.}  Configuration varieties \ 
{\bf 2.1}  Definitions \
{\bf 2.2}  Desingularizaton \ \
{\bf 3.}  The Case of $GL(n)$ \
{\bf 3.1}  Subset families \
{\bf 3.2}  Chamber families \
{\bf 3.3}  Varieties and defining equations \ \
{\bf 4.}  Schur and Weyl modules \
{\bf 4.1}  Definitions \
{\bf 4.2}  Borel-Weil theory \
{\bf 4.3}  Demazure's character formula  \ \
{\bf 5.}  Schubert polynomials \ \
{\bf 6.}  Appendix:  Non-reduced words \ \ \ \
References
}

\section{Bott-Samelson varieties}

\subsection{Three constructions}
\label{Three constructions}

In this section, $G$ is a reductive algebraic 
group.  Our constructions are all valid over
an arbitrary field, or over the integers, but
we will use the complex numbers $\CC$ for convenience.

Let $W$ denote the Weyl group 
generated by simple reflections
$s_1, \ldots , s_r$, where $r$ is the rank of $G$.
For $w \in W$, \, $\ell(w)$ denotes the length
of a reduced (i\.e\. minimal)
decompostion $w = s_{i_1} \ldots s_{i_l}$,
and $w_0$ is the element of maximal length.

We let $B$ be a Borel subgroup, $T \subset B$ a 
maximal torus (Cartan subgroup), 
and $U_{\al} \subset B$ the one-dimensional 
unipotent subgroup associated to the root $\al$.
Let $P_k \supset B$ be the {\em minimal} parabolic
associated to the simple reflection $s_k$,
so that $P_i/B \cong \PP^1$, the projective line.
Also, take $\Ph_k \supset B$ to be 
the {\em maximal} parabolic
associated to the reflections 
$s_1,\ldots, \widehat{s_k},\ldots, s_r$. 
Finally, we have the Schubert variety 
as a $B$-orbit closure 
inside the flag variety:
$$
X_w = \overline{BwB} \subset G/B
$$  

For what follows, we fix a reduced decompostion 
of some $w \in W$,
$$
w = s_{i_1} \ldots s_{i_l},
$$ 
and we denote $\ii = (i_1,\ldots,i_l)$.

Now let $P \supset B$ 
be any parabolic subgroup of $G$,
and $X$ any space with $B$-action.  Then the 
{\em induced $P$-space} is the quotient
$$
P \mtimes^B X \eqdef (P \times X)/B
$$
where the quotient is by
the free action of $B$ on $P \times X$ given by 
$(p,x) \cdot b = (pb, b^{-1}x)$.
(Thus $(pb,x) = (p,bx)$ in the quotient.)
The key property of this construction is that
$$
\begin{array}{ccc}
X & \rightarrow & P \mtimes^B X \\
& & \downarrow \\
& & P/B
\end{array}
$$
is a fiber bundle with fiber $X$ and base $P/B$.
We can iterate this construction for a sequence
of parabolics $P, P',\ldots$,
$$
P\mtimes^B P' \mtimes^B \cdots \eqdef
P\mtimes^B (\, P' \mtimes^B (\cdots)\,).
$$

Then the {\bf quotient Bott-Samelson variety} 
of the reduced word $\ii$ is
$$
\Qii \eqdef P_{i_1} \mtimes^B \cdots \mtimes^B P_{i_l} /B.
$$
Because of the fiber-bundle property of induction,
$\Qii$ is clearly a smooth, irreducible 
variety of dimension $l$.  It is a subvariety of
$$
X_l \eqdef 
\underbrace{G\mtimes^B 
\cdots \mtimes^B G}_{l\ \mbox{\footnotesize  factors}} /B.
$$
$B$ acts on these spaces by multiplying
the first coordinate:
$$
b \cdot (p_1,p_2,\ldots,p_l) 
\eqdef (bp_1, p_2, \ldots, p_l).
$$

The original purpose of the Bott-Samelson
variety was to desingularize the Schubert
variety $X_w$ via the multiplication map:
$$
\begin{array}{ccc}
\Qii & \rightarrow & X_w \subset G/B \\
(p_1,\ldots,p_l) & \mapsto & p_1 p_2 \cdots p_l B,
\end{array}
$$
a birational morphism.

Next, consider the {\em fiber product}
$$
G/B \mtimes_{G/P} G/B
\eqdef \{(g_1,g_2) \in (G/B)^2 \mid g_1 P = g_2 P \}.
$$
We may define the 
{\bf fiber product Bott-Samelson variety}
$$
\Fii \eqdef eB \mtimes_{G/P_{i_1}} G/B 
\mtimes_{G/P_{i_2}} \cdots
\mtimes_{G/P_{i_l}} G/B \subset (G/B)^{l+1}.
$$
We let $B$ act diagonally on $(G/B)^{l+1}$;
that is, simultaneously on each factor:
$$
b \cdot (g_0 B,g_1 B,\ldots,g_l B) 
\eqdef (b g_0 B, b g_1 B, \ldots, b g_l B).
$$
This action restricts to $\Fii$.
The natural map to the flag variety is
the projection to the last coordinate:
$$
\begin{array}{ccc}
\Fii & \rightarrow & G/B \\
(e B, g_1 B, \ldots, g_l B) & \mapsto & g_l B
\end{array}
$$

Finally, let us 
define the {\bf $B$-orbit Bott-Samelson
variety} as the closure 
(in either the Zariski or analytic topologies)
of the orbit of a point $\zii$:
$$
\Oii \eqdef \overline{B \cdot \zii}
\subset G/\Ph_{i_1} \times \cdots \times G/\Ph_{i_l},
$$
where  
$$
\zii = (s_{i_1} \Ph_{i_1}, \, s_{i_1}\! s_{i_2} \Ph_{i_2} 
\, , \ldots , \, s_{i_1}\!\! \cdots\!  s_{i_l} \Ph_{i_l})
$$
Again, $B$ acts diagonally.
In this case the map to $G/B$ is more difficult to
describe, but see Sec. 
\ref{Varieties and defining equations}.

\subsection{Isomorphism theorem}
\label{Isomorphism theorem}

The three types of Bott-Samelson
variety are isomorphic.
\begin{thm}
\label{isomorphism theorem}
(i)  Let
$$
\begin{array}{cccc}
\phi : & X_l & \rightarrow & (G/B)^{l+1} \\
      & (g_1, g_2, \ldots , g_l) & \mapsto &
        (\overline{e}, \overline{g_1}, \, \overline{g_1 g_2} \, , \ldots , \, 
           \overline{g_1 g_2\! \cdots\! g_l}) ,
\end{array}
$$
where $\overline{g}$ means the coset of $g$.
Then $\phi$ restricts to an isomorphism of $B$-varieties
$$
\phi : \Qii \stackrel{\sim}{\rightarrow} \Fii. 
$$ 
(ii) Let 
$$
\begin{array}{cccccccccc}
\psi : & X_l & \rightarrow &  
            G/\Ph_{i_1} & \times & G/\Ph_{i_2} 
            & \times & \cdots & \times & G/\Ph_{i_l} \\
      & (g_0,g_1, \ldots , g_l) & \mapsto &
        (\ \ \overline{g_1}&,& \overline{g_1 g_2} &,& \ldots &,&
           \overline{g_1 g_2\! \cdots\! g_l}) ,
\end{array}
$$
where $\overline{g}$ means the coset of $g$.
Then $\psi$ restricts to an isomorphism of $B$-varieties
$$
\psi : \Qii \stackrel{\sim}{\rightarrow} \Oii.
$$
\end{thm}
{\em Proof.}  (i)  It is trivial to verify
that $\phi$ is a $B$-equivariant isomorphism from 
$X_l$ to $eB \times (G/B)^l$ and that
$\phi(\Qii) \subset \Fii$, 
so it suffices to show the reverse inclusion.
Suppose 
$$
z_f = (eB, g_1 B, \ldots, g_l B) \in \Fii .
$$
Then 
$$
z_q = \phi^{-1}(z_f) =
(g_1, g_1^{-1} g_2, g_2^{-1} g_3, \ldots ) \in X_l.
$$
By definition, $e P_{i_1} = g_1 P_{i_1}$,
so $g_1 \in P_{i_1}$.  
Also $g_1 P_{i_2} = g_2 P_{i_2}$, so
$g_1^{-1} g_2 \in P_{i_2}$, and similarly
$g_{k-1}^{-1} g_k \in P_{i_k}$.  Hence $z_q \in \Qii$,
and $\phi(z_q) = z_f$. \\
(ii) First let us show that 
$\psi$ is injective on $\Qii$.
Suppose 
$\psi(p_1, \ldots, p_l) = 
\psi(q_1, \ldots, q_l)$ for $p_k, q_k \in P_{i_k}$.
Then $p_1 \Ph_{i_1} = q_1 \Ph_{i_1}$,
so that $p_1^{-1} q_1 \in \Ph_{i_1} \cap P_{i_1} = B$.
Thus $q_1 = p_1 b_1$ for $b_1 \in B$.
Next, we have 
$$
p_1 p_2 \Ph_{i_2} = q_1 q_2 \Ph_{i_2}
= p_1 b_1 q_2 \Ph_{i_2},
$$
so that $p_2^{-1} b_1 q_2 \in \Ph_{i_2} \cap P_{i_2} = B$,
and 
$q_2 = b_1^{-1} p_2 b_2$ for $b_2 \in B$.
Continuing in this way, we find that 
\begin{eqnarray*}
(q_1, q_2, \ldots , q_l) 
& = & (p_1 b_1, b_1^{-1} p_2 b_2, 
\ldots, b_{l-1}^{-1} p_l b_l) \\
& = & (p_1, p_2, \ldots, p_l) \in X_l
\end{eqnarray*}
Thus $ \psi$ is injective on $\Qii$.

Since we are working with algebraic morphisms,
we must also check that $ \psi$ is 
injective on tangent vectors of $\Qii$.
Now, the degeneracy locus 
$$
\{z \in \Qii \mid 
\mbox{\rm Ker}\ d\psi_z \neq 0 \}
$$
is a $B$-invariant, closed subvariety of $\Qii$, and
by Borel's Fixed Point Theorem it must contain
a $B$-fixed point.  But it is easily seen that
the degenerate point 
$$
z_0 = (e,\ldots,e) \in X_l
$$ 
is the only fixed point
of $\Qii$.  Thus if $d \psi$ is injective
at $z_0$, then the degeneracy locus is empty, 
and $d \psi$ is injective on each tangent space.
The injectivity at $z_0$ is easily shown 
by an argument completely analogous to that
for global injectivity given above, 
but written additively in terms of
Lie algebras instead of
multiplicatively with Lie groups.

Thus it remains to show surjectivity:
that $\psi$ takes $\Qii$ onto $\Oii$.  
Consider
$$
\ziiq = (s_{i_1},\ldots, s_{i_l}) \in X_l,
$$ 
a well-defined point in $\Qii$.  Then 
$$
\psi(\ziiq) = \zii =
(s_{i_1} \Ph_{i_1}, s_{i_1} s_{i_2} \Ph_{i_2}, \ldots),
$$
and $\psi$ is $B$-equivariant,
so that $\psi(\Qii) \supset \
\psi (\overline{B \cdot \ziiq}) = 
\overline{B \cdot \zii} = \Oii$.

Now we need only show that 
$\psi(\Qii) \subset \Oii$, 
which results from the following:
\begin{lem}
$B \cdot \ziiq$ is an open dense orbit in $\Qii$.
\end{lem}
{\bf Proof.}
Since $\Qii$ is irreducible of dimension $l$, it
suffices to show that the orbit has (at least)
the same dimension.
We may see this 
by determining $\Stab_B(\ziiq)$.
Suppose
$$
(b s_{i_1},\ldots,s_{i_l}) 
= (s_{i_1}b_1,\, b_1^{-1}s_{i_2} b_2,
\ldots,b_{l-1}^{-1} s_{i_l} b_l) \in \Qii.
$$
Then $s_{i_l} = b_{l-1}^{-1} s_{i_l} b_{l}$,
and $b_{l-1} \in B \cap s_{i_l} B s_{i_l}$.
Repeating this calculation leftward,
we find that $b \in B \cap w B w^{-1}$,
so that $\Stab_B(\zii) \subset B \cap w B w^{-1}$.
(Recall $w = s_{i_1} \dots s_{i_l}$.)
Thus, using some well-known facts
(see \cite{Spr}) we have:
\begin{eqnarray*}
\dim( B \cdot \ziiq ) & = & \dim(B) - \dim(\Stab_B(\zii)) \\
             & \geq & \dim(B) - \dim(B \cap wBw^{-1}) \\
             & = & \dim(B)-(\dim(B) - \ell(w)\,) \\
             & = & \ell(w) \, = \, l.
\end{eqnarray*}
Since the orbit can have dimension no bigger than $l$,
we must have equality.
Thus the Lemma and the Theorem both follow.
$\bullet$ 
\begin{cor}
For $w = s_{i_1} \cdots s_{i_l}$, we have
$$
\Stab_B(\zii \in \Zii) = \Stab_B(wB \in G/B) 
= B \cap w B w^{-1}.
$$
\end{cor}

\subsection{Open cells}
\label{Open cells}

In view of the Theorem, we will let
$\Zii$ denote the abstract Bott-Samelson
variety defined by any of our three versions.
It contains the degenerate $B$-fixed point $z_0$
defined by:
\begin{eqnarray*}
z_0 & = & (e,e,\ldots) \in \Qii \\
&=& (eB,eB,\ldots) \in \Fii \\
&=& (e \Ph_{i_1},e \Ph_{i_2},\ldots) \in \Oii
\end{eqnarray*}
as well as the generating $T$-fixed point
whose $B$-orbit is dense in $\Zii$:
\begin{eqnarray*}
\zii &=& (s_{i_1}, s_{i_2}, s_{i_3}, \ldots) \in \Qii \\
&=& (e B, s_{i_1}  B, s_{i_1} s_{i_2} B, \ldots) 
\in \Fii \\
&=& (s_{i_1} \Ph_{i_1}, s_{i_1} s_{i_2} \Ph_{i_2},
\ldots) \in \Oii
\end{eqnarray*}

We may parametrize the dense 
orbit $B \cdot \zii \subset \Zii$ by an affine cell.
Consider the normal ordering of the positive
roots associated to the reduced word $\ii$.
That is, let 
$$
\be_1 = \al_{i_1},\ \be_2 = s_{i_1}( \al_{i_2}),\
\be_3 =  s_{i_1} s_{i_2} (\al_{i_3}),\ \cdots
$$
Recall that $U_{\be_k}$ is the one-dimensional unipotent
subgroup of B corresponding to the positive root $\be_k$.
Then we have a direct product:
$$
B = U_{\be_1}\! \cdots U_{\be_l}\cdot (B\cap wBw^{-1}),
$$
so that the multiplication map 
$$
\begin{array}{ccc}
U_{\be_1} \times \cdots \times U_{\be_l} 
& \rightarrow & B\cdot \zii \\
(u_1,\ldots,u_l) & \mapsto & u_1\cdots u_l \cdot \zii
\end{array}
$$
is injective, and an isomorphism of varieties.
The left-hand side is isomorphic to an affine
space $\CC^l$. 
\\[1em]
$\Zii$ also contains an opposite big cell 
centered at $z_0$ which is not the orbit of a group.
Consider the 
one-dimensional unipotent
subgroups $U_{-\al_{i}}$ corresponding to the 
negative simple roots $-\al_i$.  
The map
$$
\begin{array}{cccc}
\CC^l\ \cong \!\! & U_{-\al_{i_1}}\! \times \cdots \times U_{-\al_{i_l}}
& \rightarrow & \Qii \\
& (u_1,\ldots,u_l) & \mapsto & 
(u_1 , \dots , u_l)
\end{array}
$$
is an open embedding.

In the case of $G = GL(n)$, $B = $ upper triangular
matrices, we may write an element of
$U_{-\al_{i_k}}$ as $u_k = I+t_k e_k$,
where $I$ is the identity matrix, $e_k$ is
the sub-diagonal
coordinate matrix $e_{(i_k+1,i_k)}$,
and $t_k \in \CC$.
If we further map $\Qii$ to $G/B$ via the
natural multiplication map, we get
$$
\begin{array}{ccc}
(t_1,\ldots,t_l) & \mapsto & 
(I+t_1 e_1) \cdots (I+t_l e_l) \\
\CC^l & \rightarrow & N_- \\
\cap && \cap \\
\Qii &\rightarrow & G/B \\
(p_1,\ldots,p_l) & \mapsto & p_1\cdots p_l B
\end{array}
$$
where $N_-$ denotes the unipotent lower triangular
matrices (mod $B$).  Thus the multiplication on the
bottom is a compactification of the matrix factorizations
studied by Berenstein, Fomin, and Zelevinsky \cite{BFZ}.

\section{Configuration varieties}

We define a class of varieties (more general
than the Schubert varieties) which are 
desingularized by Bott-Samelson varieties. 

\subsection{Definitions}
\label{Definitions}

We continue with 
the case of a general reductive group $G$.  
Given a sequence of 
Weyl group elements $\ww = (w_1,\ldots,w_k)$ 
and a sequence
of indices $\jj = ({j_1},\ldots,{j_k})$,
we consider the $T$-fixed point 
$$
z_{\ww \jj} = (w_1 \Ph_{j_1},\ldots,w_k \Ph_{j_k})
\in G/ \Ph_{j_1} \times \cdots \times G/ \Ph_{j_k},
$$
and we define the {\em configuration variety}
as the $G$-orbit closure
$$
\FF_{\ww \jj} \eqdef 
\overline{ G \cdot z_{\ww \jj}}
\subset G/ \Ph_{j_1} \times \cdots \times G/ \Ph_{j_k}.
$$
$G$ acts on this variety by multiplying each factor
simultaneously (the diagonal action).

We may define a ``flagged'' version of this 
construction by replacing $G$ with $B$.
The {\em flagged configuration variety}
is the $B$-orbit closure
$$
\FFB_{\ww \jj} \eqdef 
\overline{ B \cdot z_{\ww \jj}}
\subset G/ \Ph_{j_1} \times \cdots \times G/ \Ph_{j_k}.
$$
Again, $B$ acts diagonally.
\\[1em]
{\bf Examples.} (a) Take 
$\ww = (w,w,\ldots,w)$ 
for any $w \in W$ and
$\jj = (1,2,\ldots,r)$ \, (where $r = \rank G$).
Then the configuration variety is isomorphic to
the flag variety of $G$, and 
the flagged configuration variety is isomorphic
to the Schubert variety of $w$: 
$$ 
\FF_{\ww \jj} \cong G/B \ \ \ \ \ \ \ \ \ \ \ 
\FFB_{\ww \jj} \cong X_w \ .
$$
(b) For $\jj = \ii = (i_1,i_2,\ldots)$, 
a reduced word, and 
$\ww = (s_{i_1},  s_{i_1} s_{i_2},\ldots)$,
the flagged configuration variety is exactly our
orbit version of the Bott-Samelson variety: 
$\FFB_{\ww \jj} = \Oii = \Zii.$ \ 
$\bullet$
\\[1em]
{\bf Remark.} For a given $G$, there are 
only finitely many configuration varieties up to
isomorphism.  
In fact, suppose a list ($\ww$, $\jj$) has
repetitions of some element of $\ww$ with 
identical corresponding entries in $\jj$.
Then we may remove the repetitions
and the configuration variety will not 
change (up to $G$-equivariant isomorphism), 
only the embedding.
Thus, all configuration varieties
are projections of a maximal variety.  
This holds for the flagged and unflagged cases.
\\[1em]
{\bf Example.}  The maximal configuration 
variety for $G = GL(3)$ is the {\em space of
triangles} \cite{MTri}, and corresponds to
$$
\begin{array}{ccccccccl}
\ww & = & (e,& e,& s_1,& s_2,& s_2 s_1,& s_1 s_2 &) \\
\jj & = & (1,& 2,& 1,  & 2,  &    1,   &    2 &).
\end{array}
$$
Further entries would be redundant: for example,
$s_1 \Ph_2 = e \Ph_2$. 
All other configuration varieties are
obtained by omitting some entries of $\ww$ 
and the corresponding entries
of $\jj$.  Hence there are at most
$2^6$ configuration varieties for $G$.
$\bullet$
\\[1em]
One might attempt to broaden the definition
of configuration varieties
by replacing the minimal homogeneous spaces
$G/\Ph_j$ by $G/P$ for arbitrary
parabolics $P \supset B$.   
This gives the same class of varieties,
however, since any $G/P$ can be embedded
equivariantly inside a product of $G/\Ph_j$'s,
resulting in isomorphic orbit closures.
Once again, this changes only the embeddings,
not the varieties.

Varieties similar to our $\FF_{\ww \jj}$ are defined 
and some small cases are analyzed 
in Langlands' paper \cite{Langlands}.

\subsection{Desingularization}
\label{Desingularization}

Very little is known about
general configuration varieties.
However, certain of them are
well understood because they can be 
desingularized by Bott-Samelson varieties.

Recall that a sequence 
$\ww = (w_1,\ldots,w_K)$ of Weyl group elements
is increasing in the weak order on $W$
if there exist $u_1,u_2,\ldots, u_K$ such that 
$w_k = u_1 u_2 \cdots u_k$ and
$\ell(w_k) = \ell(w_{k-1}) + \ell(u_k)$ 
for all $k$. 

For $\ww = (w_1,\ldots,w_K)$ and $\jj = (j_1,\ldots,j_K)$,
let $\ww^+ = (e,\ldots,e,w_1,\ldots,w_k)$ with $r$  added
entries of $e$, and $\jj^+ = (1,2,\ldots,r,j_1,\ldots,j_K)$.
Clearly
$$
\FFB_{\ww,\jj} \cong \FFB_{\ww^+ \jj^+}.
$$
\begin{prop}
\label{desingularization}
If $\ww$ is increasing in the weak order
and $\jj$ is arbitrary,
then the flagged configuration variety
$\FFB_{\ww \jj}$ can be desingularized
by a Bott-Samelson variety.  That is, 
there exists a reduced word $\ii$ 
and a regular birational morphism 
$$
\pi: \Zii \rightarrow \FFB_{\ww \jj}.
$$

Furthermore, the unflagged configuration
variety $\FF_{\ww^+\jj^+}$
is desingularized by the composite map
$$
G\mtimes^B \Zii \stackrel{\mbox{id} \times \pi}{\rightarrow}
G\mtimes^B \FFB_{\ww \jj} \cong 
G\mtimes^B \FFB_{\ww^+ \jj^+} \stackrel{\mu}{\rightarrow}
\FF_{\ww^+ \jj^+},
$$
where $\mbox{id} \times \pi$ is the map induced from $\pi$,
and $\mu$ is the multiplication map
$(g, v) \mapsto g \cdot v$.
\end{prop}
{\bf Remark.} The map 
$$
G\mtimes^B \Zii \rightarrow G\mtimes^B \FFB_{\ww \jj}
\rightarrow \FF_{\ww\jj}
$$
is a surjection from a smooth space to $\FF_{\ww\jj}$,
but it is not birational in general.
We will see in Sec 4 that for the purposes of
Borel-Weil theory, this map can substitute for a
desingularization of $\FF_{\ww\jj}$.
$\bullet$
\\[1em]
To prove the Proposition, we will need the following
\begin{lem}  (a) For any $w \in W$ and parabolic $P$
with Weyl group $W(P)$,
we have a unique factorization 
$w = \tw y$, where $y \in W(P)$,
$\tw$ has minimal length in $\tw W(P)$,
and $\ell(w) = \ell(\tw)+\ell(y)$.
\\
(b) Suppose $w \in W$ has minimum length in the coset
$w W(P)$, and consider the points 
$wP \in G/P$ and $wB \in G/B$.
Then $\Stab_B(wP) = \Stab_B(wB)$.
\end{lem}
{\bf Proof of Lemma.}  
(a) Well-known (see \cite{Humphreys}, \cite{Hiller}). \\
(b) The $\supset$ containment is clear, so
we prove the other.
Let $\Del$ denote the set of roots of $G$,
$\Del_+$ the positive roots, $\Del(P)$ the roots of $P$,
etc.  From considering the corresponding Lie algebras
we obtain: 
$$
\begin{array}{rcl}
\dim \Stab_{B}(wB) & = & | \Del_+ \cap w(\Del_+) | \\
\dim \Stab_{B}(wP) & = & 
| \Del_+ \cap w(\Del_+ \!\cup \Del(P)) | .
\end{array}
$$
But the two sets on the right
are identical.  In fact, if $w$ is minimal in $wP$,
then $\Del_+ \cap w(\Del_-(P)) = \emptyset$.  (See
\cite{Humphreys}, 5.5, 5.7.)
$\bullet$.

\mbox{}\\[1em]
{\bf Proof of Proposition.}
Denote $W_k \eqdef  W(\Ph_{j_k})$, a
parabolic subgroup of the Weyl group.
Given $\ww$ and $\jj$, 
we define a new sequence
$\tww = (\tw_1,\ldots,\tw_K)$.
Take $\tw_k$ to be 
the minimum-length coset representative
in $w_k W_k$, so that 
$w_k = \tw_k y_k$ for some $y_k \in W_k$.
I claim the new sequence $\tww$ is still
increasing in the weak order.  In fact,
if $w_k = u_1\cdots u_k$ and 
$\tu_k$ is minimal in $u_k W_k$, then
$\tw_k = \tw_{k-1} y_k \tu_k$ 
and $\ell(\tw_k) = \ell(\tw_{k-1}) + \ell(y_k)
+ \ell(\tu_k)$.
Note that it is possible 
that $\tu_k = e$, and 
$\tw_k = \tw_{k+1}$.

Now let $\ii$ be any reduced decomposition
of the increasing sequence $\tww$:  that is,
for each $k$ we have a reduced decompostion
$\tw_k = s_{i_1} s_{i_2} \cdots s_{i_{l(k)}}$, 
where $l(k) = \ell(\tw_k)$,
so that 
$0 \leq l(1) \leq l(2) \cdots \leq l(K)=l$.
Also, $i_{l(k)} = j_k$ for all $k$.  
%It follows that $\tw_k = \tw_{k-1}$, or
%the last letter $i_p$ in any such word must be
%$s_{j_k}$ (otherwise $\tw_k$ could not be minimal
%in $\tw_k W_k$).
Define a projection map from the Bott-Samelson variety
to the configuration variety:
$$
\begin{array}{cccc}
\phi: &\Zii = \Oii & \rightarrow & \FFB_{\ww \jj} \\
& (g_1 \Ph_{i_1},\ldots,g_l \Ph_{i_l}) & \mapsto & 
(g_{l(1)} \Ph_{j_1} \, ,
\ldots, \,   g_{l(K)} \Ph_{j_K}).
\end{array}
$$
I claim $\phi$ is well-defined, $B$-equivariant,
onto, regular, and birational.  

Now, $\tw_k$ and $w_k$ are equal modulo $W_k$,
so $\tw_k \Ph_{j_k} = w_k \Ph_{j_k}$, 
and thus
$$
\phi(z_{\ii}) = z_{\tww \jj} = z_{\ww \jj} 
\in \FFB_{\ww \jj}.
$$
Since $\Zii = \overline{B \cdot z_{\ii}}$, this 
implies that the image of $\phi$ lies inside
$\FF_{\ww \jj}$, and $\phi$ is well-defined.
It is clearly $B$-equivariant and therefore onto
(since $\FFB_{\ww \jj}$ is a $B$-orbit closure).

The map is regular, and to show it is birational we
need only check that it is a bijection between the big
$B$-orbits in the domain and image.  That is, we must
show equality of the stabilizers
$$
\Stab_B(\zii) = \Stab_B(z_{\tww \jj}).
$$ 
By the corollary in Section \ref{Isomorphism theorem}, 
we have $\Stab_B(\zii)=\Stab_B(wB \in G/B)$
for $w=s_{i_1}\cdots s_{i_l(K)}=\tw_K$.

Now we use induction on
the length of the sequence $\ww$.
If the length $K=1$, we have immediately
that $\Stab_B(z_{\ww \jj}) = 
\Stab_B(\tw_K \Ph_{j_K}) = \Stab_B(\tw_K B)$
by the above Lemma.
Assuming the assertion for $\ww'=(w_1,\ldots,w_{K-1})$
and using the Lemma, we have 
$$
\begin{array}{rcl}
\Stab_B(z_{\tww \jj}) & = &
\Stab_B(z_{\tww' \jj}) \cap \Stab_B(\tw_K \Ph_{j_K}) \\
&=& \Stab_B(\tw_{K-1}B) \cap \Stab_B(\tw_K B) \\
&=& \Stab_B(\tw_K B) .
\end{array}
$$
\mbox{} \\[1em]
The remaining assertions about
the unflagged $\FF_{\ww^+ \jj^+}$ follow easily.
That is, the map of fiber bundles
$$
G \mtimes^B \Zii
\rightarrow G\mtimes^B \FFB_{\ww^+\jj^+}
$$
is $G$-equivariant, onto, and
regular and birational by our results above, 
and so is the multiplication map
$$
G \mtimes^B  \FFB_{\ww^+ \jj^+} 
\rightarrow \FF_{\ww^+ \jj^+}
$$
since $\Stab_G(z_{\ww^+\jj^+}) = \Stab_B(z_{\ww^+\jj^+})$.
$\bullet$

\section{The Case of $GL(n)$}

We begin again, restating many of our results 
more explicitly for the general linear group
$G = GL(n,\CC)$.  In this case $B$ = upper triangular
matrices, $T =$ diagonal matrices, $r = n-1$,
$$
P_k = \{ (x_{ij}) \in GL(n) \mid x_{ij} = 0 
\mbox{ if } i>j \mbox{ and } (i,j) \neq (k+1,k) \},
$$
$$
\Ph_k = \{ (x_{ij}) \in GL(n) \mid x_{ij} = 0 
\mbox{ if } i>k\geq j \},
$$
and $G/\Ph_k \cong \Gr(k,\CC^n)$, the
Grassmannian of $k$-dimensional subspaces
of complex $n$-space.

Also $W$ = permutation matrices, $\ell(w) =$
the number of inversions of a permutation $w$,
$s_i =$ the transposition $(i,i+1)$, and
the longest permutation is $w_0 = n, n-1, \ldots, 2,1$.  
We will frequently use the notation
$$
[k] = \{1,2,3,\ldots,k\}.
$$

\subsection{Subset families}
\label{Subset families}

First, we introduce some combinatorics.
Define a {\em subset family} to
be a collection $D = \{C_1,C_2,\ldots\}$
of subsets $C_k \subset [n]$.
The order of the subsets is irrelevant
in the family, and we do not allow subsets
to be repeated.

This relates to the previous sections as follows.
To a list of permutations $\ww = (w_1,\ldots,w_K)$,
$w_k \in W$, and a list of indices 
$\jj = (j_1, \ldots, j_K)$, $1 \leq j_k \leq n$,
we associate a subset family:
$$
D = D_{\ww \jj} \eqdef 
\{ w_1 [j_1], \ldots, w_K[j_K] \} .
$$
Here $w[j] = \{ w(1), w(2), \ldots, w(j)\}$.

Now suppose the list of indices 
$\ii = (i_1, i_2, \ldots , i_l)$ encodes 
a reduced decomposition 
$w = s_{i_1} s_{i_2} \cdots s_{i_l}$ of
a permutation into a minimal number of
simple transpositions.
We let $\ww = (s_{i_1}, s_{i_1} s_{i_2}, \ldots, w)$
and $\jj = \ii$, and we define the 
{\em reduced chamber family} 
$D_{\ii} \eqdef D_{\ww \jj}$.

Further, define the {\em full chamber family}
$$
\Dtii \eqdef \{ [1], [2], \ldots, [n] \} \cup \Dii ,
$$
(which is $D_{\ww^+ \jj^+}$ in our previous notation).

We tentatively connect these structures
with geometry. 
Let $\CC^n$ have the standard basis
$e_1,\ldots,e_n$.
For any subset $C = \{j_1,\ldots, j_k\} \subset [n]$,
the coordinate subspace 
$$
E_C = \Span_{\CC}\{e_{j_1},\ldots,e_{j_k}\}\,  
\in\, \Gr(k)
$$
is a $T$-fixed point in a Grassmannian.
A subset family 
corresponds to a $T$-fixed point
in a product of Grassmannians 
$$
z_D = (E_{C_1},E_{C_2},\ldots)
\in \Gr(D) \eqdef \Gr(\,|C_1|\,) \times 
\Gr(\,|C_2|\,) \times \ldots.
$$
This is consistent with our  previous notation for an
arbitrary $G$:  for $D = D_{\ww \jj}$,
we have $z_D = z_{\ww \jj}$.
We defined configuration varieties and
Bott-Samelson varieties as orbit
closures of such points (see also below, 
Sec \ref{Varieties and defining equations}).
\\[1em]
{\bf Examples.} 
For $n = 3$, $G = GL(3)$, $\ii = \jj = 121$,
we have $\ww = (s_1,s_1\! s_2,\, s_1\! s_2 s_1)$,
and the reduced chamber family
\begin{eqnarray*}
D_{121} & = & 
\{\, s_1[1],\, 
     s_1 s_2[2],\, 
         s_1 s_2 s_1[1]\, \}\\
&=&  \{\,\{2\},\{2,3\},\{3\}\,\} \\
&=& \{2,23,3\} 
\end{eqnarray*}
The full chamber family is
$D^+_{121} =  \{1,12,123,2,23,3\}$.
The chamber family of the other reduced word 
$\ii = 212$ is $D_{212} = \{13, 3, 23\}$.
%\\[.5em]
%For $n=4$, the word $\ii = 312132$ is a reduced
%decomposition of the longest permutation $w_0$,
%and we have 
%$\Dii = \{124,2,24,4,234,34\}$.
\\[.5em]
For $n=4$, let $\ww = (e,s_1,s_1,s_3 s_2,s_1)$,
$\jj = (2,1,3,1,1)$.  Then we have the
subset family
$$
\begin{array}{ccl}
\Dwj &=& \{ \,
e[2],\, s_1[1],\, s_1[3],\, 
s_3 s_2[1],\, s_1[1]
\,\} \\
&=& \{12,2,123,3,2\} = \{12,123,2,3\}
\end{array}
$$
Note that we remove repetitions in $D$.
The associated $T$-fixed configuration is
$$
z_D = (E_{12},E_{123},E_{2},E_{3})
\in \Gr(D) = 
\Gr(2) \times \Gr(3) \times \Gr(1) \times \Gr(1).
$$
$\bullet$

\subsection{Chamber families}
\label{Chamber families}

Chamber families have a rich structure.
(See \cite{LZ}, \cite{MFour}, \cite{RSNew}.)
Given a full chamber family $\Dtii$, we may 
omit some of its elements 
to get a subfamily $D \subset \Dtii$.
The resulting {\em chamber subfamilies}
can be characterized as follows.  

For two sets $S, S' \subset [n]$, 
we say $S$ is {\em elementwise less than} $S'$,
$S \lelt S'$, 
if $s < s'$ for all
$s \in S$, $s' \in S'$.
Now, a pair of subsets $C, C' \subset [n]$
is {\em strongly separtated} if
$$
(C \setminus C') \lelt (C' \setminus C)
\ \ \ \mbox{or} \ \ \ 
(C' \setminus C) \lelt (C \setminus C') \ ,
$$
where $C\setminus C'$ denotes 
the complement of $C'$ in $C$.
A family of subsets is called strongly separated
if each pair of subsets in it is strongly separated.

\begin{prop}{(LeClerc-Zelevinsky \cite{LZ})\ }
\label{LeClerc-Zelevinsky}
A family $D$ of subsets of $[n]$ is
a chamber subfamily, $D \subset \Dtii$
for some $\ii$,
if and only if $D$ is
strongly separated.
\end{prop}
{\bf Remarks.} (a) Reiner and Shimozono \cite{RSNew}
give an equivalent description of strongly separated
families.  Place the subsets of the family into
lexicographic order.  
Then $D = (C_1 \leqlex C_2 \leqlex \cdots)$ is 
strongly separated if and only if it is 
``\%-avoiding'':  that is, if 
$i_1 \in C_{j_1}$, 
$i_2 \in C_{j_2}$ with $i_1 > i_2$, $j_1 < j_2$,
then $i_1 \in C_{j_2}$ or $i_2 \in C_{j_1}$.  \\
(b) If $\ii = (i_1,\ldots,i_l)$ is
an initial subword of $\ii' = (i_1, \ldots, i_l, 
\ldots, i_N)$, then $\Dii \subset D_{\ii'}$.
Thus the chamber families associated to
decompositions of the longest permutation $w_0$
are the maximal strongly separated families.\\
(c) In \cite{MFour}, we describe the ``orthodontia''
algorithm to determine
a reduced decomposition $\ii$ associated to a given
strongly separated family.  See also \cite{RSNew}.
\\[1em]
{\bf Examples.} (a)   
For $n=3$, the chamber families
$D^+_{121} = \{1,12,123,2,23,3\}$ and 
$D^+_{212} = \{1,12,123,13,3,23\}$
are the only maximal strongly separated families.  
The sets $13$ and $2$ are the only pair
not strongly separated from each other. 
\\
(b) For $n=4$,  the strongly separted family
$D=\{24,34,4\}$ is contained in the
chamber sets of the reduced words
$\ii = 312132$ and $\ii=123212$.
$\bullet$
\\[1em]
Chamber families can be represented pictorially
in several ways, one of the most natural being 
due to Berenstein, Fomin, and Zelevinsky \cite{BFZ}.
The {\em wiring diagram} or {\em braid diagram}
of the permutation
$w$ with respect to the reduced word $\ii$
is best defined via an example.

Let $G = GL(4)$, $w = w_0$ 
(the longest permutation),
and $\ii = 312132$.
On the left and right ends 
of the wiring diagram are
the points 1,2,3,4 in two columns.
Each point $i$ on the left is connected to
the point $w(i)$ on the right by a curve 
which is horizontal and disjoint from the other
curves except for certain crossings.  
The crossings, 
read left to right, correspond 
to the entries of $\ii$. 
The first entry 
$i_1 = 3$ corresponds to 
a crossing of the curve on level 3 with the
one on level 4. (The other curves continue
horizontally.)  The second entry $i_2 = 1$
crosses the curves on level 1 and 2, and so on.
\\[1em]
\centerline{\Large \bf FIGURE 1 }
\\[1em]
If we add crossings only up to the 
$l^{\mbox{\tiny th}}$ step,
we obtain the wiring diagram of the truncated 
word $s_{i_1} s_{i_2} \cdots s_{i_l}$.

Now we may construct the chamber family 
$$
\Dtii = (1,12,123,1234,124,2,24,4,234,34)
$$
as follows.  Label each of the curves of the
wiring diagram by its point of origin on the left.
Into each of the connected 
regions between the curves,
write the numbers of those curves
which pass above the region.  
Then the sets of numbers inscribed in these chambers
are the members of the family $\Dtii$. 
If we list the chambers from left to right, 
we recover the natural order in which these subsets 
appear in $\Dtii$. 

Another way to picture 
a chamber family, or any subset family, 
is as follows.  
We may consider a subset 
$C = \{j_1,j_2,\ldots\} \subset [n]$
as a column of $k$ squares in the rows $j_1, j_2,\ldots$.
For each subset $C_k$ in the chamber family, 
form the column associated to it, and
place these columns next to 
each other.  The result is an array
of squares in the plane called a 
{\em generalized Young diagram}.

For our word $\ii = 312132$, we draw the
(reduced) chamber family as:
$$
\Dii \ = \ 
%(124,2,24,4,234,34) \ = \ 
\begin{array}{ccccccc}
1 & \Box &      &      &      &      &      \\
2 & \Box & \Box & \Box &      & \Box &      \\   
3 &      &      &      &      & \Box & \Box \\
4 & \Box &      & \Box & \Box & \Box & \Box
\end{array}
$$
where the numbers on the left of the diagram
indicate the level.
See \cite{RS1}, \cite{MNW}, \cite{MFour}.

\subsection{Varieties and defining equations}
\label{Varieties and defining equations}

To any subset family $D$ we have associated
a $T$-fixed point in a product of Grassmannians,
$z_D \in \Gr(D)$, and we may define as before
the {\em configuration variety} of $D$
to be the closure of the $G$-orbit of $z_D$:
$$
\FFD = \overline{G\cdot z_D} \subset \Gr(D);
$$
and the {\em flagged configuration variety}
to be the closure of its $B$-orbit:
$$
\FFBD = \overline{B\cdot z_D} \subset \Gr(D).
$$
Furthermore, if $D = \Dii$, a chamber family,
then the {\em Bott-Samelson variety} is the flagged
configuration variety of $\Dii$:
$$
\Zii = \Oii = \FFB_{\Dii}.
$$
(We could also use the full chamber family
$\Dtii$, since the extra coordinates correspond
to the standard flag fixed under the $B$-action.)

Thus $\FFD$, $\FFBD$, and $\Zii$ 
can be considered as
varieties of configurations of 
subspaces in $\CC^n$, 
like the flag and Schubert varieties.
We will give defining equations for 
the Bott-Samelson varieties analogous to 
those for Schubert varieties.

For a family $D$, 
define the {\em flagged inclusion variety}
$$
\II^B_D = \left\{
\begin{array}{c|c}
(V_C)_{C\in D}\in \Gr(D) &
\begin{array}{c} \forall\, C, C' \in D,\ \,
C \subset C' \Rightarrow V_C \subset V_{C'} \\
\mbox{and  }\ \forall\, [i] \in D, \ \ V_{[i]} = \CC^i
\end{array}
\end{array}
\right\}.
$$
$B$ acts diagonally on $\II^B_D$.
\\[1em] 
{\bf Example.} 
For $n=4$, $\ii=312132$, we may use the picture in the
above example to write the inclusion variety 
$\II^B_{\Dii^+}$ as
the set of all 10-tuples of subspaces of $\CC^4$
$$
(V_1,V_{12},V_{123},V_{1234},V_{124},V_2,V_{24},V_4,
V_{234},V_{34})
$$
with $\dim( V_C\!) = |C|$ and satisfying the following
inclusions:
$$
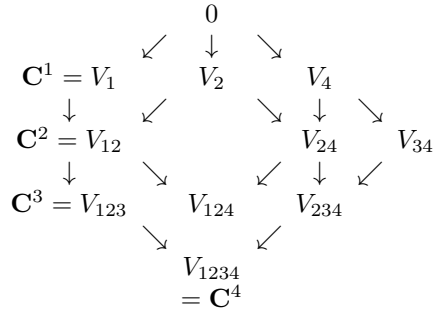
\begin{array}{ccccccc}
&&0&&&&\\
&\!\!\!\!\!\!\swarrow\!\!\!\!\!\!&\!\downarrow\!&\!\!\!\!\!\!\searrow\!\!\!\!\!\!&&&\\
\CC^1=V_1&&V_2&&V_4&&\\
\!\downarrow\!&\!\!\!\!\!\!\swarrow\!\!\!\!\!\!&&\!\!\!\!\!\!\searrow\!\!\!\!\!\!&\downarrow&\!\!\!\!\!\!\searrow\!\!\!\!\!\!&\\
\CC^2=V_{12}&&&&V_{24}&&V_{34}\\
\!\downarrow\!&\!\!\!\!\!\!\searrow\!\!\!\!\!\!&&\!\!\!\!\!\!\swarrow\!\!\!\!\!\!&\downarrow&\!\!\!\!\!\!\swarrow\!\!\!\!\!\!\\
\CC^3=V_{123}&&V_{124}&&V_{234}&&\\
&\!\!\!\!\!\!\searrow\!\!\!\!\!\!&&\!\!\!\!\!\!\swarrow\!\!\!\!\!\!&&&\\
&&V_{1234}&&&& \\
&& = \CC^4 &&&&
\end{array}
$$
where the arrows indicate inclusion of subspaces.
\begin{thm}
For every reduced word $\ii$, we have
$\Zii \cong \II^B_{\Dii^+}$.
\end{thm}
{\bf Proof.}
Note that the generating point
$z_{\Dii^+}$ lies in $\II^B_{\Dii^+}$, 
and $\II^B_{\Dii^+}$ is
B-equivariant, so $\Zii \subset \II^B_{\Dii^+}$.

To show the reverse inclusion, we use
our previous characterization 
$$
\Zii \cong \Fii = e \mtimes_{G/P_{i_1}}
G/B \mtimes_{G/P_{i_2}} G/B \mtimes_{G/P_{i_3}}
\cdots   \mtimes_{G/P_{i_l}} G/B.
$$
We may write this variety as the
$(l+1)$-tuples of flags
$ (V_1^{(k)} \subset V_2^{(k)} \subset 
\cdots \subset \CC^n) $, \, $k = 0,1,\ldots,l$,
such that: $V_i^{(k)} = V_i^{(k+1)}$ for all
$k$ and all $i \neq i_k$; 
and $V_i^{(0)} = \CC^i$ for all $i$. 

Consider the map
$$
\begin{array}{cccc}
\theta: & \Fii & \rightarrow & \Gr(D) \\[.1em]
& (V_1^{(k)} \subset V_2^{(k)} \subset \cdots)_{k=0}^l
& \mapsto & (V_{i_1}^{(1)}, V_{i_2}^{(2)}, \ldots)
\end{array}
$$
We have seen in Theorem \ref{isomorphism theorem} that
$\Zii = \Oii = \IM(\theta)$, 
since $\theta = \psi \circ \phi^{-1}$.  
It remains to show that 
$\II^B_{\Dii^+} \subset \IM(\theta)$.

For each $k$, define 
$k^- = \max\{ m \mid m < k,\, i_m = i_k +1 \}$
and $k^+ = \min\{ m \mid m > k,\, i_m = i_k +1 \}$.
Then it is easily seen that a configuration
$(V_1,V_2,\ldots) \in \Gr(D)$ lies in
$\IM(\\theta)$ exactly when: \\
(i) for each $k$, we have $V_k \subset V_{k^-}$
and $V_k \subset V_{k^+}$
provided $k^-$ or $k^+$ is defined; \\
(ii) for each $k$, if $k^-$ is not defined, 
then $V_k \subset \CC^{i_k+1}$; and \\
(iii) for each $i$, if $k = \min\{m \mid i_m = i+1\}$,
then $\CC^i \subset V_k$.

Note that for any $k$, the $k^{\mbox{\tiny th}}$ subset
of $\Dii$ is 
$$
\begin{array}{rcl}
C_k & = & s_{i_1} \cdots s_{i_k} [i_k] \\
& = & s_{i_1} \cdots s_{i_k} 
\cdots s_{i_{k^+}} [i_k] \\
& \subset &  s_{i_1} \cdots s_{i_k} 
\cdots s_{i_{k^+}} [i_k+1] \\
& = &  s_{i_1} \cdots s_{i_k} 
\cdots s_{i_{k^+}} [i_{k^+}] \\
&=& C_{k^+}
\end{array}
$$
We can write similar inclusions of subsets for
the other conditions (i)-(iii).
This shows that the inclusions defining $\II^B_{\Dii^+}$
do indeed imply those defining $\IM(\theta)$,
Q.E.D.
$\bullet$

\begin{conj}
For any subset family $D$,
a configuration $(V_C)_{C\in D} \in \Gr(D)$
lies in $\FF_D$ exactly if, for every subfamily
$D' \subset D$, 
$$
\dim(\bigcap_{C \in D'} V_C) \geq \, |\! \cap_{C \in D'} C |
$$
$$
\dim(\sum_{C \in D'} V_C) \leq \, |\! \cup_{C \in D'} C |
$$
\end{conj}
Note that a configuration 
$(V_1,\ldots,V_l) \in \Gr(D)$ lies
in the flagged configuration variety
$\FFB_D$ if and only if 
$(\CC^1,\ldots, \CC^n, V_1,\ldots, V_l)$
lies in the unflagged variety $\FF_{D^+}$
of the augmented diagram $D^+ \eqdef \{[1],
[2], \ldots [n]\} \cup D$.  Hence the
above conjecture gives conditions defining
flagged configuration varieties as well as unflagged.
\\[1em]
{\bf Examples.}
(a) If $D=\Dii$ is a chamber family, the conjecture
reduces to the previous Theorem. 
\\
(b) The conjecture is known  if
$D$ satisfies the ``northwest condition''
(see \cite{MNW}):
that is, the elements of $D$ can be arranged
in an order $C_1, C_2, \ldots$ such that
if $i_1 \in C_{j_1}$, 
$i_2 \in C_{j_2}$, 
then $\min(i_1,i_2) \in C_{\min(j_1,j_2)}$.
In fact, it suffices in this case to consider only
the  intersection conditions of the conjecture.
$\bullet$
\\[1em]
It would be interesting to know whether the 
determinantal equations implied by the conditions
of the above Theorem and Conjecture define
$\FFD \subset \Gr(D)$ scheme-theoretically.

Now, let $D$ be a strongly separated family.
We know by Proposition \ref{LeClerc-Zelevinsky} that $D$
is part of some chamber family $\Dii$, 
and by Theorem
\ref{desingularization} we may take $\ii$
so that the projection map $\Zii = \FFB_{\Dii}
\rightarrow \FFB_D$ is birational.
\\[1em]
{\bf Example.}
Let $n=7$, and consider the family $D$ consisting
of the single subset $C = 12457$.  Its configuration
variety is the Grassmannian $\FFD = \Gr(5, \CC^7)$,
and its flagged configuration variety is the Schubert
variety 
$$
\FFB_D = X_{211} = \{V \in \Gr(5) \mid
\CC^2 \subset V,\, \dim(\CC^5 \cap V) \geq 4 \}.
$$
By the orthodontia algorithm \cite{MFour}, we find that
this is desingularized by the reduced
word $\ii = 3465$, for which 
$\Dii = \{124,1245,123457, 12457\}$ and
$$
\Zii = 
\left\{
\begin{array}{c}
(V_{124},V_{1245},V_{123457},V_{12457}) 
\in \Gr(3) \times \Gr(4)
\times \Gr(6) \times \Gr(5)  \\[.1cm]

\mbox{  such that  }\ \ \ 
\CC^2 \subset V_{124} \subset \CC^4 \subset V_{123457}\ \ , \ \
V_{1245} \subset \CC^5\ , \\
\mbox{} \ \ \ \ \ \ \
V_{124} \subset V_{1245} \subset V_{12457} \subset V_{123457}
\end{array}
\right\}.
$$
The desingularization map is the projection
$$
\pi:(V_{124},V_{1245},V_{123457},V_{12457}) 
\mapsto V_{12457} .
$$
In \cite{MNW} and Zelevinsky's work \cite{Zel},
there are given several other desingularizations of
Schubert varieties, all of them expressible 
as configuration varieties.
$\bullet$

\section{Schur and Weyl modules}

We relate generalized Schur and 
Weyl modules for $GL(n)$,
which are defined in completely elementary terms,  
to the sections of line bundles on configuration
varieties, and hence to the coordinate rings of these
varieties.

One the one hand, 
this yields an unexpected Demazure character
formula for the Schur modules, including the
skew Schur functions and Schubert polynomials. 
On the other hand, it gives an
elementary construction for line-bundle
sections on Bott-Samelson varieties.

\subsection{Definitions}
\label{Schur and Weyl Modules: Definitions}

We have associated to any subset family
$D = \{C_1,\ldots,C_k\}$ a configuration 
variety $\FF_D$ with $G$-action, and a
flagged configuration variety $\FFB_D$
with $B$-action.  
Now, assign an integer multiplicity
$\mm(C) \geq 0$ to each subset $C \in D$.
For each pair $(D,\mm)$, we define a $G$-module
and a $B$-module, which will turn out to sections
of a line bundle on $\FF_D$ and $\FFB_D$.

In the spirit of DeRuyts 
\cite{F} and Desarmenien-Kung-Rota \cite{DKR},
we construct these ``Weyl modules'' $M_{D,\mm}$
inside the coordinate ring of $n\times n$ matrices,
and their flagged versions $\MB_{D,\mm}$
inside the coordinate ring of upper-triangular matrices.  
(I am grateful to Mark Shimozono for pointing out this
form of the definition.)

Let $\CC[x_{ij}]$  
(resp.  $\CC[x_{ij}]_{i\leq j} $ )
denote the polynomial
functions in the variables $x_{ij}$ with
$i,j \in [n]$\
(resp. $x_{ij}$ with $1 \leq i\leq j\leq n$).
For $R, C \subset [n]$ with $|R|=|C|$,
let 
$$
\Del_C^R = \det(x_{ij})_{(i\in R, j\in C)} \in \CC[x_{ij}]
$$
be the minor determinant of the matrix $x = (x_{ij})$
on the rows $R$ and the columns $C$.
Further, let 
$$
\tDel_C^R = \Del_C^R |_{x_{ij} = 0 , \ \forall\, i>j} 
\in \CC[x_{ij}]_{i\leq j}
$$
be the same minor evaluated on an upper triangular
matrix of variables.

Now, for a subset family
$D=\{C_1,\dots,C_l\}$, $\mm=(m_1,\ldots,m_l)$,
define the {\em Weyl module}
$$
M_{D,\mm} = \Span_{\CC}\left\{\
\Del_{C_1}^{R_{11}}\cdots \Del_{C_1}^{R_{1m_1}}
 \Del_{C_2}^{R_{21}} \ldots \Del_{C_l}^{R_{lm_l}}
\left| 
\begin{array}{c}
\forall\, k,\! m \ \
R_{km} \subset [n] \\[.1cm]
\mbox{ and }\ |R_{km}|=|C_k| 
\end{array}
\right.
\right\}.
$$
That is, a spanning vector is a product of
minors with column indices equal to the elements
of $D$ and row indices taken arbitrarily.

For two sets $R = \{i_1,\ldots,i_c\}$,
$C = \{j_1,\ldots,j_c\}$ we say
$R \leqcomp C$ (componentwise inequality) if
$i_1 \leq j_1$, $i_2\leq j_2$, \ldots.
Define the {\em flagged Weyl module}
%$$
%\MB_{D,\mm} = \Span_{\CC}\left\{
%\prod_{C \in D} \prod_{k = 1}^{m(C)}
%\Del_C^{R_{Ck}} 
%\ \ \left| \ \ 
%\begin{array}{c}
%R_{Ck} \subset [n], \\
%|R_{Ck}|=|C| \\
%R_{Ck} \leqcomp C
%\end{array}
%\right.
%\right\}.
%$$
$$
\MB_{D,\mm} = \Span_{\CC}\left\{\
\tDel_{C_1}^{R_{11}}\cdots \tDel_{C_1}^{R_{1m_1}}
 \tDel_{C_2}^{R_{21}} \ldots \tDel_{C_l}^{R_{lm_l}}
\left| 
\begin{array}{c}
\forall\, k,m \ \ 
R_{km} \subset [n] \\[.1cm]
 |R_{km}|=|C_k| ,\, 
R_{km} \leqcomp C_k
\end{array}
\right.
\right\}.
$$

For $f(x) \in \CC[x_{ij}]$, a matrix $g \in G$
acts by left translation, $(g \cdot f)(x) = f(g^{-1}x)$.
It is easily seen that this restricts to a $G$-action
on $M_{D,\mm}$ and similarly we get 
a $B$-action on $\MB_{D,\mm}$.

We clearly have the diagram of $B$-modules:
$$
\begin{array}{ccc}
M_{D,\mm} & \subset & \CC[x_{ij}] \\
\downarrow & & \downarrow \\
\MB_{D,\mm} & \subset & \CC[x_{ij}]_{i\leq j}
\end{array}
$$
where the vertical maps ($x_{ij} \mapsto 0$
for $i>j$) are surjective.  
That is, $\MB_{D,\mm}$ is a quotient of $M_{D,\mm}$.

The {\em Schur modules} are defined to be the duals
$$
S_{D,\mm} \eqdef (M_{D,\mm})^* \ \ \ \ \ \
S^B_{D,\mm} \eqdef (\MB_{D,\mm})^*.
$$
We will deal mostly with the Weyl modules, but
everything we say will of course 
also apply to their duals.
\\[1em]
{\bf Example.} We adopt the ``Young diagram''
method for picturing subset families.  
(See Sec \ref{Chamber families}.) 
Let $n=4$, $D = \{234,34,4\}$, $m = (2,0,3)$.
(That is, $m(234)=2$, $m(34)=0$, $m(4)=3$.)
We picture this by writing each column repeatedly,
according to its multiplicity.  Zero multiplicity
means we omit the column.  Thus
$$
(D,\mm) = 
\begin{array}{cccccc}   
1 &      &      &      &      &       \\
2 & \Box & \Box &      &      &       \\
3 & \Box & \Box &      &      &       \\
4 & \Box & \Box & \Box & \Box & \Box 
\end{array}
\ \ \ \ 
\tau = 
\begin{array}{c|ccccc}   
1 &      &      &      &      &       \\
2 & 1 & 1 &      &      &       \\
3 & 3 & 2 &      &      &       \\
4 & 4 & 3 & 2 & 4 & 3 
\end{array}
$$
The spanning vectors for $M_{D,\mm}$
correspond to all column-strict fillings of this
diagram by indices in $[n]$. 
For example, the filling $\tau$ above corresponds to
$$
\Del_{234}^{134}\
\Del_{234}^{123}\
\Del_{4}^{2}\
\Del_{4}^{4}\
\Del_{4}^{3}\
$$
$$
= \left| \begin{array}{ccc}
x_{12} & x_{13} & x_{14} \\
x_{32} & x_{33} & x_{34} \\
x_{42} & x_{43} & x_{44} 
\end{array} \right| \cdot
\left| \begin{array}{ccc}
x_{12} & x_{13} & x_{14} \\
x_{22} & x_{23} & x_{24} \\
x_{32} & x_{33} & x_{34} 
\end{array} \right| \cdot
x_{24}\cdot x_{44}\cdot x_{34}
$$
$$
= \left( \begin{array}{ccccc|ccccc}
 1 & 1 &   &   &   & 2 & 2 &   &   &   \\
 3 & 2 &   &   &   & 3 & 3 &   &   &   \\
 4 & 3 & 2 & 4 & 3 & 4 & 4 & 4 & 4 & 4
\end{array} \right)
$$
The last expression is in the letter-place 
notation of Rota et al \cite{DKR}.

A basis may be extracted from this spanning
set by considering only the row-decreasing fillings
(a normalization of the semi-standard tableaux), 
and in fact the
Weyl module is the dual of the classical Schur 
module $S_{\lam}$ associated to the shape $D$
considered as the Young diagram $\lam = (5,2,2,0)$.

The spanning elements of the flagged Weyl 
module $\MB_{D,\mm}$ correspond to the
``flagged'' fillings of the diagram: those for
which the number $i$ does not appear above the
$i^{\mbox{\tiny th}}$ level.   
For the diagram above, all the
column-strict fillings are flagged, and
$M_{D,\mm} \cong \MB_{D,\mm}$.  

However, for 
$$
(D',\mm) = 
\begin{array}{cccccc}   
1 &      &      &      &      &       \\
2 & \Box & \Box &      &      &       \\
3 & \Box & \Box & \Box & \Box & \Box  \\
4 & \Box & \Box &      &      &
\end{array}
$$
$$
\tau_1 = 
\begin{array}{c|ccccc}   
1 &   &   &   &   &   \\
2 & 2 & 1 &   &   &   \\
3 & 3 & 2 & 4 & 3 & 4 \\
4 & 4 & 3 &   &   &   
\end{array}
\ \ \ \ 
\tau_2 = 
\begin{array}{c|ccccc}   
1 &   &   &   &   &   \\
2 & 2 & 1 &   &   &   \\
3 & 3 & 2 & 3 & 2 & 3 \\
4 & 4 & 4 &   &   &   
\end{array}
$$
the filling $\tau_1$ is {\em not} flagged,
since 4 appears on the 3rd level, but $\tau_2$
{\em is} flagged, and corresponds to the
spanning element
$$
\tDel_{234}^{234}\
\tDel_{234}^{124}\
\tDel_{3}^{3}\
\tDel_{3}^{2}\
\tDel_{3}^{3}
= \left| \begin{array}{ccc}
x_{22} & x_{23} & x_{24} \\
0 & x_{33} & x_{34} \\
0 & 0 & x_{44} 
\end{array} \right| \cdot
\left| \begin{array}{ccc}
x_{12} & x_{13} & x_{14} \\
x_{22} & x_{23} & x_{24} \\
0 & 0 & x_{44} 
\end{array} \right| \cdot
x_{33}\cdot x_{23}\cdot x_{33}.
$$

We have $M_{D,\mm} \cong M_{D',\mm} \cong \MB_{D,\mm}
\cong S^*_{(5,2,2,0)}$,
the dual of a classical (irreducible) Schur module
for $GL(4)$,
and $\MB_{D',\mm} \cong S^*_{(2,5,2,0)}$,
the dual of the Demazure module 
with lowest weight $(0,2,5,2)$
and highest weight $(5,2,2,0)$.
Cf. \cite{RS1}, \cite{MFour}.
$\bullet$ 
\\[1em]
{\bf Remarks.} 
(a) In \cite{LM} we make a general definition
of ``standard tableaux'' giving bases of
the Weyl modules for strongly separated
families. \\
(b)  We briefly indicate the equivalence 
between our definition of the Weyl modules
and the tensor product definition given in 
\cite{ABW}, \cite{RS1}, \cite{MNW}.

Let $Y = Y_{D,\mm} \subset \NN \times \NN$ be
the generalized Young diagram of squares in the
plane associated to $(D,\mm)$ as in the above
examples, and let $U = (\CC^n)^*$.
One defines $\Mt_Y = U^{\otimes Y} \gamma_Y$,
where $\gamma_Y$ is a generalized Young symmetrizer.
The spanning vectors $\Del_{\tau}$ of $M_{D,\mm}$
correspond to the fillings $\tau:Y \rightarrow [n]$.
Then the map 
$$
\begin{array}{ccc}
M_{D,\mm} & \rightarrow & \Mt_{D,\mm} \\
\Del_{\tau} & \mapsto &
\left( \bigotimes_{(i,j)\in Y} e^*_{\tau(i,j)} \right)
\gamma_Y
\end{array}
$$
is a well-defined isomorphism of $G$-modules, and
similarly for the flagged versions.
This is easily seen from the definitions,
and also follows from the Borel-Weil theorems proved
below and in \cite{MNW}.

\subsection{Borel-Weil theory}
\label{Borel-Weil theory}

A configuration variety $\FF_D \subset \Gr(D)$
has a natural family of line bundles defined
by restricting the determinant or Plucker bundles
on the factors of $\Gr(D)$.
For $D = (C_1,C_2,\ldots)$, and multiplicities
$\mm = (m_1, m_2, \ldots)$, we define
$$
\begin{array}{ccc}
\LLmm & \subset & \OO(m_1,m_2,\ldots) \\
\downarrow & & \downarrow \\
\FF_D & \subset & \Gr(D) = \Gr(|C_1|) \times \Gr(|C_2|)
\times \cdots
\end{array}
$$
We denote by the same symbol $\LLmm$ this line
bundle restricted to $\FFB_D$.
Note that in the case of a Bott-Samelson variety
$\FF_D = \Zii$, this is the well-known line bundle
$$
\LLmm \cong 
{P_{i_1} \times \cdots \times P_{i_l} \times \CC 
\over B^l}
$$
$$
(p_1,\ldots,p_l,v)\cdot (b_1,\ldots,b_l)
\eqdef
(p_1 b_1,\ldots,b_{l-1}^{-1}p_l b_l, \,
\om_{i_1}(b_1^{-1})^{m_1} \cdots \om_{i_l}(b_l^{-1})^{m_l}\,v),
$$
$\om_i$ denoting the fundamental weight 
$\om_i(\diag(x_1,\ldots,x_n)) =
x_1 x_2 \cdots x_i$.

Note that if $m_k \geq 0$ for all $k$ 
(resp. $m_k >0$ for all $k$) then
$\LLmm$ is effective (resp. very ample).
However, $\LLmm$ may be effective even if
some $m_k < 0$.  See \cite{LM}.

\begin{prop}
Let $(D,\mm)$ be a strongly separated subset family
with multiplicity.  Then we have \\
(i) $M_{D,\mm} \cong H^0(\FF_D,\LLmm)$ \\[.1em]
and $H^i(\FF_D,\LLmm) = 0$ for $i>0$. \\[.2em]
(ii) $\MB_{D,\mm} \cong H^0(\FFB_D,\LLmm)$ \\[.1em]
and $H^i(\FFB_D,\LLmm) = 0$ for $i>0$. \\[.2em]
(iii) $\FF_D$ and $\FFB_D$ are normal varieties,
projectively normal with respect to $\LLmm$, and
have rational singularities.
\end{prop}
{\bf Proof.}
First, recall that we can identify the
sections of a bundle over a single Grassmannian,
$\OO(1) \rightarrow \Gr(i)$, with 
linear combinations of
minors in the homogeneous Stiefel coordinates
$$
x = \left(\begin{array}{ccc}
x_{11} & \cdots & x_{1i} \\
\vdots & \ddots & \vdots \\
x_{n1} & \cdots & x_{ni} \\
\end{array} \right)
\in \Gr(i),
$$
namely the $i \times i$ minors $\Del^R(x)$
on the rows $R \subset [n]$, $|R| = i$.
Thus, a typical spanning element of 
$H^0(\Gr(D),\OO(\mm))$ is the section
$$
\Del^{R_{11}}(x^{(1)})
 \cdots  \Del^{R_{11}}(x^{(1)})\
\Del^{R_{21}}(x^{(2)})
\cdots  \Del^{R_{lm_l}}(x^{(l)}),
$$
where $x^{(k)}$ represents the homogeneous coordinates 
on each factor $\Gr(|C_k|)$ of $\Gr(D)$,
and $R_{km}$ are arbitrary subsets with $|R_{km}|=i_k$.

Now, restrict this section to 
$\FF_D \subset \Gr(D)$ and then further to the
dense $G$-orbit $G\cdot z_D \subset \FF_D$.
Parametrizing the orbit by $g \rightarrow g\cdot z_D$,
we pull back the resulting 
sections of $H^0(\FF_D,\LLmm)$
to certain functions on 
$G \subset \mbox{Mat}_{n\times n}(\CC)$, 
which are precisely the products
of minors defining the spanning 
set of $M_{D,\mm}$.
This shows that
$$
M_{D,\mm} \cong
\IM\left[ 
H^0(\Gr(D),\OO(\mm)) \
\stackrel{\mbox{\footnotesize rest}}{\rightarrow}
H^0(\FF_D,\LLmm)
\right].
$$
Similarly for $B$-orbits, we have
$$
\MB_{D,\mm} \cong
\IM\left[ 
H^0(\Gr(D),\OO(\mm)) \
\stackrel{\mbox{\footnotesize rest}}{\rightarrow}
H^0(\FFB_D,\LLmm)
\right].
$$

Now we invoke the key vanishing result,
\cite{MNW} Prop. 28
(due to W. van der Kallen and S.P. Inamdar, 
based on the work of O. Mathieu \cite{Mat}, P. Polo, et.al.)
The conditions $(\alpha)$ and $(\beta)$
of that Proposition apply to $\FF_D$ because
$D$ is contained in a chamber family $D^+_{\ii}$
(Prop. \ref{LeClerc-Zelevinsky} above).
Furthermore, the proof of \cite{MNW}, Prop. 28
goes through identically with $\FFB_D$ in place
of $\FF_D$, merely replacing $\FF_{w_0;u_1,\ldots,u_r}$
by $\FF_{e;u_1,\ldots,u_r}$.  

All of the assertions of our Proposition now
follow immediately from the corresponding parts
of \cite{MNW}, Prop. 28.
$\bullet$.

\begin{prop}
\label{extend by zero}
Suppose $(D,\mm)$, $(\tD, \tmm)$ are
strongly separated subset families with
$D \subset \tD$, $\tmm(C) = \mm(C)$ for $C\in D$,
$\mm(C) = 0$ otherwise.
Then the natural projection
$\pi: \Gr(\tD) \rightarrow \Gr(D)$
restricts to a surjection
$\pi: \FF_{\tD} \rightarrow \FF_D$,
and induces an isomorphism
$$
\pi^*:  
H^0(\FF_D,\LLmm)
\stackrel{\sim}{\rightarrow}
H^0(\FF_{\tD},\LL_{\tmm}),
$$
and similarly for the flagged case.
\end{prop}
{\bf Proof.}  For the unflagged case, this follows
immediately from \cite{MNW}, Prop. 28.  Again,
the argument given there goes through for the flagged
case as well. $\bullet$
\\[1em]
{\bf Remarks.}
(a) Note that the proposition holds even if 
$\dim \FF_{\tD} > \dim \FF_D$. \\
(b) The Proposition allows us to reduce
Weyl modules for strongly separated
families to those for {\em maximal} 
strongly separated families, that is
chamber families. $\bullet$
\\[1em]
We may conjecture that the results of
this section hold not only in the 
strongly separated case, but for all
subset families and configuration varieties.

\subsection{Demazure's character formula}
\label{Demazure's character formula}

We now examine how the iterative structure of
Bott-Samelson varieties influences
the associated Weyl modules.

Define Demazure's isobaric divided difference
operator
$\Lam_i : \CC[x_1,\ldots,x_n] 
\rightarrow \CC[x_1,\ldots,x_n] $,
$$
\Lam_i f = {x_i f - x_{i+1} s_i f \over x_i - x_{i+1}}.
$$
For example for $f(x_1,x_2,x_3) = x_1^2 x_2^2 x_3$,
$$
\begin{array}{rcl}
\Lam_2 f(x_1,x_2,x_3) &=&
{ x_2(x_1^2 x_2^2 x_3) - x_3(x_1^2 x_3^2 x_2) 
\over x_2 - x_3 } \\
&=& x_1^2 x_2 x_3 (x_2+x_3).
\end{array}
$$
For any permutation with a reduced decompostion
$w = s_{i_1}\ldots s_{i_l}$, define
$$
\Lam_{w} \eqdef \Lam_{i_1} \cdots \Lam_{i_l},
$$
which is known to be independent of the reduced
decomposition chosen.

By the (dual) character of a $G$- or $B$-module $M$,
we mean 
$$
\Char M = \tr(\diag(x_1,\ldots,x_n)|M^*)\,
\in\, \CC[x_1^{\pm 1},\ldots,x_n^{\pm 1}].
$$
(We must take duals to get polynomial functions
as characters.)
Let $\om_i$ denote the $i$th fundamental weight,
the multiplicative character
of $B$ defined by $\om_i(\diag(x_1,\ldots,x_n)) =
x_1 x_2 \cdots x_i$.

\begin{prop}
\label{character formula}
Suppose $(D,\mm)$ is strongly separated, and
$$
D \subset \Dii^+ = 
\{ [1],\ldots,[n],C_1,\ldots,C_l \},
$$
for some reduced word $\ii = (i_1,\ldots,i_l)$.
Define $\tmm = (k_1,\ldots,k_n, m_1,\ldots, m_l)$
by 
$\tmm(C) = \mm(C)$ for $C \in D$, $\tmm(C) = 0$ 
otherwise.
Then 
$$
\Char \MB_{D,\mm} = 
\om_1^{k_1} \cdots \om_n^{k_n}
\Lam_{i_1} \om_{i_1}^{m_1} \cdots
\Lam_{i_l} \om_{i_l}^{m_l}.
$$
Furthermore,
$$
\Char M_{D,\mm} = 
\Lam_{w_0} \Char \MB_{D,\mm} ,
$$
where $w_0$ denotes the longest permutation.
\end{prop}
{\bf Remark.}  We explain in \cite{LS} how
one can recursively generate
the standard tableaux for $\MB_D$ (in \cite{LM}) 
by ``quantizing'' this character formula. 
 See also \cite{MFour}.
\\[1em]
We devote the rest of this section to proving the Proposition.

For a subset $C = \{j_1,j_2,\ldots\} \subset [n]$, 
and a permutation $w$, 
let $wC = \{w(j_1), w(j_2),\ldots\}$,
and for a subset family $D = \{C_1,C_2,\ldots\}$,
let $wD = \{wC_1, wC_2,\ldots\}$.
Now, for $i \in [n-1]$, let 
$$
 \Lam_i D  \eqdef \{s_i[i]\} \cup s_i D  ,
$$
where $s_i [i] =  \{1,2,\ldots,i-1,i+1\}$.
We say that $D$ is {\em $i$-free} for $i\in [n]$
if for every $C\in D$, we have $C \cap \{i,i+1\}
\neq \{i+1\}$.

\begin{lem}
\label{i-free}
Suppose $(D,\mm)$ is strongly separated and $i$-free.\\
(i) $\FFB_{\Lam_i D} \cong P_i \mtimes^B \FFB_D$ .\\
(ii) $\FFB_{s_i D} \cong 
P_i \cdot \FFB_D \subset \Gr(D)$ .\\
(iii) The projection 
$\FFB_{\Lam_i D} \rightarrow \FFB_{s_i D}$
is regular, surjective, and birational. \\
(iv) Let $\tmm$ be the multiplicity on $\Lam_i D$ defined
by $\tmm(s_i C) = \mm(C)$ for $C\in D$, 
$\tmm(s_i[i]) = m_0$.
The bundle $\LL_{\tmm} \rightarrow \FFB_{\Lam_i D}$
is isomorphic to 
$$
\LL_{\tmm} \cong P_i \mtimes^B 
\left( (\om_i^{m_0})^* \otimes \LLmm \right),
$$
where $(\om_i^{m_0})^* \otimes \LLmm$
indicates the bundle $\LLmm \rightarrow \FFB_D$
with its $B$-action twisted by the multiplicative
character $(\om_i^{m_0})^* = \om_i^{-{m_0}}$.

\end{lem}
{\bf Proof.}
(i)  Since $D$ is $i$-free, we have
$U_i z_D = z_D$, where $U_i$ is the
one-dimensional unipotent subgroup corresponding
to the simple root $\al_i$.  We may factor
$B$ into a direct product of subgroups,
$B = U_i B' = B' U_i$.  Then
$$
\FFB_D = \overline{B \cdot z_D} 
= \overline{B'\cdot z_D}.
$$
Hence the $T$-fixed point 
$(s_i,z_D) \in P_i\mtimes^B \FFB_D$
has a dense $B$-orbit:
$$
\begin{array}{rcl}
\overline{B \cdot (s_i,z_D)} &=& 
\overline{(U_i B' s_i, z_D)} \\[.1em]
&=& (\overline{U_i s_i}, \overline{B' \cdot z_D}) \\[.1em]
&=& P_i \mtimes^B \FFB_D.
\end{array}
$$
Clearly, the injective map 
$$
\begin{array}{cccc}
\psi: & P_i \times^B \Gr(D) & \rightarrow & 
\Gr(i) \times \Gr(D) \\
& (p,V) & \mapsto & (p \CC^i, pV)
\end{array}
$$
takes $\psi(s_i,z_D) = z_{\Lam_i D}$,
the $B$-generating point of $\FFB_{\Lam_i D}$.
Thus $\psi:P_i \times^B \FFB_D  \\
\rightarrow \FFB_{\Lam_i D}$ is an isomorphism. \\
(ii+iii)  By the above, the projection is a bijection
on the open B-orbit, and hence is birational.
The image of the projection is $P_i \cdot \FFB_D$,
which must be closed since $P_i \mtimes^B \FFB_D$
is a proper (i.e. compact variety).  \\
(iv) Clear from the definitions.
$\bullet$

\begin{lem}
Let $(D,\mm)$ be a strongly separated family
and $i \in [n-1]$.  
Let 
$$
\begin{array}{rcl}
\FF' & = & P_i \mtimes^B \FFB_D \\
\LL' & = & P_i \mtimes^B \LLmm .
\end{array}
$$
so that $\LL' \rightarrow \FF'$ is a line
bundle.
Then
$$
\Char H^0(\FF',\LL') = \Lam_i \Char H^0(\FFB_D,\LLmm).
$$
\end{lem}
{\bf Proof.}  By Demazure's analysis of induction
to $P_i$ (see \cite{Dem1},
``construction \'{e}l\'{e}mentaire'')
we have
$$
\Lam_i \Char H^0(\FFB_D,\LLmm) = 
\Char H^0(\FF',\LL') - 
\Char H^1(\, P_i/B, H^1(\FFB_D,\LLmm)\,).
$$
However, we know by \cite{MNW}, Prop.28
that $H^0(\FFB_D,\LLmm)$ has a good filtration,
so that the $H^1$ term above is zero.
$\bullet$

\begin{cor}
If $(D,\mm)$ is strongly separated and $i$-free,
and $(\Lam_i D, \tmm)$ is a diagram with
multiplicities
$\tmm(s_i C) = m(C)$ for $C\in D$, 
$\tmm(s_i[i]) = m_0$, then 
$$
\Char \MB_{\Lam_i D,\tmm} = 
\Lam_i \om_i^{m_0} \Char M^B_{D,\mm}.
$$
If $m_0 = 0$, then
$$
\Char \MB_{s_i D,\mm} 
= \Char \MB_{\Lam_i D,\tmm}
= \Lam_i \Char M^B_{D,\mm}
$$
\end{cor}
This follows immediately from the above
Lemmas and Proposition \ref{extend by zero}.
\mbox{}\\[1em]
{\bf Proof of Proposition.} 
The first formula of the 
Proposition now follows from the
above Lemmas and Prop \ref{extend by zero}.
The second statement follows from Demazure's 
character formula, combined with the vanishing
result of \cite{MNW} Prop.28. $\bullet$.
%\\ AAAAAAAAAAAAAAAAAAAAAAAAAAAAAAAAAAA \\
%Amplify!

\section{Schubert polynomials}

In this section, we again work with $G =GL(n)$.
As a general reference, see Fulton \cite{F}.
\\[.1em]

There are two classical computations of 
the singular cohomology ring $H^.(G/B, \CC)$
of the flag variety.  That of Borel 
identifies the cohomology with a coinvariant
algebra 
$$
c: H^.(G/B, \CC) \stackrel{\sim}{\rightarrow}
 \CC[x_1,\ldots,x_n]/I_+,
$$
where $I_+$ is the the ideal generated by
the non-constant symmetric polynomials.
The map $c$ is an isomorphism of 
graded $\CC$-algebras, and the generator
$x_i$ represents the Chern class of the 
$i^{\mbox{\tiny th}}$ quotient of the 
tautological flag bundle.
(This is not the dual of an effective divisor.)

The alternative picture of Schubert gives as a linear
basis for $H^.(G/B, \CC)$ the Schubert classes
$\sigma_w = [X_{w_0 w}]$, the Poincare duals
of the Schubert varieties.

The isomorphism between these pictures
was defined by Bernstein-Gelfand-Gelfand
\cite{BGG} and by Demazure \cite{Dem1},
and given a precise combinatorial form by
Lascoux and Schutzenberger \cite{LS}.
It identifies certain {\em Schubert polynomials}
$\SS(w) \in \CC[x_1,\ldots,x_n]$ with
$c(\sigma_w) = \SS(w)\, (\mbox{mod } I_+)$,
and enjoying many remarkable properties.

They can be defined combinatorially by
a descending recurrence, starting with
the representative of the fundamental class
of $G/B$.  For any permutation $w$
with $ws_i < w$ in the Bruhat order, and
$w_0$ the longest permutation, we have
$$
\SS(w_0) = x_1^{n-1} x_2^{n-2} \cdots x_{n-2}^2 x_{n-1}
$$
$$
\SS(w s_i) = \partial_i \SS(w),
$$
where we use the divided difference operator
$\partial_i: \CC[x_1,\ldots,x_n] \rightarrow
 \CC[x_1,\ldots,x_n]$,
$$
\partial_i f = {f - s_i f \over x_i - x_{i+1}}.
$$
(Note that $\Lam_i = \partial_i x_i$.
This is special to the root system of type
$A_{n-1}$.)
\\[.5em]
{\bf Example.}  For $G=GL(3)$, 
we have $\SS(w_0) = x_1^2 x_2$,
$\SS(s_1 s_2) = x_1 x_2$,
$\SS(s_2 s_1) = x_1^2$,
$\SS(s_2) = x_1 + x_2$,
$\SS(s_1) = x_1$,
$\SS(e) = 1$.
$\bullet$
\\[1em]
To compute any $\SS(w)$, we write 
$w_0 = w s_{i_1} \cdots s_{i_r}$ for some
reduced word $s_{i_1} \cdots s_{i_r}$,
and we have 
$$
\SS(w) = \partial_{i_1} \cdots \partial_{i_r} (x_1^{n-1}
x_2^{n-2} \cdots x_{n-1}).
$$
In particular, we may take $i_k$ to be the
{\em first ascent} of 
$w_k = w s_{i_1} \cdots s_{i_{k-1}}$;
that is, $i_k =$ the smallest $i$ 
such that $w_k(i+1)>w_k(i)$.

We now give a completely different geometric
interpretation of the polynomials
$\SS(w)$ in terms of configuration varieties
and Weyl modules.
For a permutation $w$ define the
{\em inversion family} 
$I(w) = \{C_1(w),\ldots,C_{n-1}(w)\}$
with
$$
C_j(w) = \{ i \in [n] \ \mid\, i<j, \ w(i)>w(j) \}
$$
We may write this in our usual form 
$(D,\mm)$ by dropping any of the $C_j(w)$ 
which are empty,
and counting identical sets 
with multiplicity.
We use the same symbol $I(w)$
to denote this multiset $(D,\mm)$, 
so that $I(w) - {C}$ means we
decrease by one the multiplicity 
of the element $C \in I(w)$.
It is well-known that $I(w)$ is strongly
separated.  (In fact, it is northwest.  See 
\cite{RS1}, \cite{RS3}, \cite{MNW})

\begin{thm}{(Kraskeiwicz-Pragacz \cite{KP})}
$$
\Char \MB_{I(w)}  = \SS(w) .
$$
\end{thm}
{\bf Proof.} (Magyar-Reiner-Shimozono)\,  
Let $\chi(w) = \Char \MB_{I(w)}$.
We must show that $\chi(w)$ satisfies the
defining relations of $\SS(w)$.

First, $I(w_0) = \{[1],\ldots,[n-1]\}$,
$$
\MB_{I(w_0)} = 
\CC \cdot \tDel_1^1 \tDel_{12}^{12} \ldots 
\tDel_{[n-1]}^{[n-1]},
$$
a one-dimensional $B$-module, and
$ \chi(w_0) = x_1^{n-1} x_2^{n-2} \cdots x_{n-1} $.

Now, suppose $ws_i < w$, and $i$ is the first
ascent of $w s_i$. 
Then the $w(i)^{\mbox{th}}$ 
element of $I(w)$ 
is $C_{w(i)}(w) = [i]$.  Letting
$$
I'(w) \eqdef I(w) - \{\, [i]\, \},
$$
it is easily seen that: \\
(i) $I'(w)$ is $i$-free, \\
(ii) $I(w) = I'(w)\, \cup \{\, [i]\, \}$, and \\
(iii) $I(w s_i) = s_i I'(w) \, \cup \, \{\,[i-1]\,\} $.\\
(Set $[0] = \emptyset$.)
%Case (i): Suppose $j \neq w(i)$.  
%Then $C_j(ws_i) = s_i C_j(w)$.  Also,
%we cannot have $i+1 \in C_j(w)$, $i \not\in C_j(w)$,
%because $i$ is a descent of $w$.  
%Thus $I'(w) \eqdef I(w) - \{C_{w(i)}(w)\}$ is $i$-free,
%and $I(ws_i) - \{C_{w(i)}(ws_i)\} = s_i I'(w)$.
%Case (ii):  Suppose $j = w(i)$.  Then
%$C_j(w s_i) = [i-1]$, $C_j(w) = [i]$.

Hence we obtain trivially:
%\\ AAAAAAAAAAAAAAAAAAAAAAAA \\
%by Lemma ??, amplify! paraphrase the previous paragraph
$$
\begin{array}{rcl}
\chi(w) &=& x_1\! \cdots x_i \, \Char \MB_{I'(w)}  \\[.2em]
\chi(ws_i) &=& 
x_1\! \cdots x_{i-1}\, \Char \MB_{s_i I'(w)}.
\end{array}
$$

Since $I'(w)$ is strongly separated and $i$-free,
Cor 14 implies that
$$
\Char \MB_{s_i I'(w)} = \Lam_i \Char \MB_{I'(w)}.
$$
This is the key step of the proof.

Thus we have
$$
\begin{array}{rcl}
\chi(ws_i) & = & (x_1 \cdots x_{i-1})\ \Lam_i \Char \MB_{I'(w)} \\
&=& \Lam_i x_i^{-1} (x_1 \cdots x_i)\, \Char\MB_{I'(w)} \\
&=& \Lam_i x_i^{-1}\, \chi(w) \\
&=& \partial_i\, \chi(w) .
\end{array}
$$

But now, using the
the first-ascent sequence to write
$w_0 = w s_{i_1} \cdots s_{i_r}$,
we compute
$$
\chi(w) = \partial_{i_1} \cdots \partial_{i_r} (x_1^{n-1}
x_2^{n-2} \cdots x_{n-1}) = \SS(w).
$$
$\bullet$

Our Demazure character formula 
(Prop \ref{character formula}) now
allows us to compute Schubert polynomials
by a completely different recursion from the
usual one.  In particular, the defining recursion
goes from higher to lower degree, whereas our
Demazure formula goes from lower to higher.
\\[1em]
{\bf Example.}
For the permutation $w = 24153$ in $GL(5)$,
we have $I(w) = \{ 12, 24 \}$  (neglecting the
empty set).  Then the first-ascent sequence 
gives us:  
$$
\SS(w) = \partial_1  \partial_3
  \partial_2  \partial_1  \partial_4 \partial_3 
(x_1^4 x_2^3 x_3^2 x_4).
$$
However, it is easier to compute
that $I(w) \subset \Dii^+$ for a chamber
family with $\ii = 132$ (= reduced word $s_1 s_3 s_2$),
so that $\Dii = \{2, 124, 24\}$ and 
$$
\begin{array}{rcccccccl}
\Dii^+ = & \{1,& 12,& 123,& 1234,& 12345,& 2,& 124 & 24\} \\
\mm = & ( 0, &    1,&   0,&   0, &    0, & 0,&  0,& 1)
\end{array}
$$
$$
\begin{array}{rcl}
\SS(w) & = & x_1 x_2\, \Lam_1\, \Lam_3\, \Lam_2\, (x_1 x_2) \\
       & = & x_1 x_2\, (x_1 x_2 + x_1 x_3 + x_1 x_4 + 
              x_2 x_3 + x_2 x_4)
\end{array}
$$
See \cite{MFour} for more examples of such computations.
$\bullet$

\section{Appendix: Non-reduced words}

Let $G$ again be an arbitrary reductive group of rank $r$.

For future reference, we note that many of our
results hold when the decomposition
$w=s_{i_1}\cdots s_{i_l}$
is not of minimal length (that is,
$\ell(w) < l$).  
We call the resulting $\ii = (i_1,\ldots,i_l)$
(with $i_k \in \{1,\ldots,r\}$)
a {\em non-reduced} word.

In this case the quotient and fiber product definitions
of the Bott-Samelson variety apply without change,
and we still have $\Zii \cong \Qii \cong \Fii$, as shown in
Thm 1(i).
However, $\Zii$ is no longer the $B$-orbit closure
of a $T$-fixed point, so we can no longer define
$\Oii$.  Nevertheless, the map
$$
\psi:X_l \rightarrow \Gr_G(\ii) \eqdef 
G/\Ph_{i_1}\times \cdots G/\Ph_{i_l} 
$$
of Thm 1(ii) is still
injective on $\Qii \subset X_l$ 
(the first part of the proof of Thm 1(ii) 
is unchanged).  Thus we may define an ``embedded''
version of $\Zii$,
$$
\Gii \eqdef \psi(\Qii) \subset \Gr_G(\ii) ,
$$
so that $\Gii = \Oii$ if $\ii$ is reduced.

We can also define analogues of
Weyl modules for a general $G$ and $\ii$.
We once again have the minimal-degree line
bundles $\OO(1)$ over the $G/\Ph_i$, and
hence $\OO(\mm) = \OO(m_1,\ldots,m_l)$ over
$\Gr_G(\ii)$.  Let $\LLmm$ be the restriction
of $\OO(\mm)$ to $\Gii$.  Then 
define 
$$
\MB_{\ii,\mm} \eqdef H^0(\Zii,\LLmm).
$$

These modules no longer embed in $\CC[B]$,
but they do have a spanning set of Pl\"ucker coordinates,
the restrictions of sections from the ambient
space $\Gr_G(\ii)$:
\begin{prop}
Let $\ii= (i_1,\ldots,i_l)$ be
an arbitrary word (not necessarily reduced),
and $\mm = (m_1,\ldots,m_l)$ with $m_j \geq 0$
for all $j$.

Then the restriction map 
$$
H^0(\Gr_G(\ii),\OO(\mm)) \rightarrow H^0(\Zii,\LLmm)
$$
is surjective.
Furthermore, $H^i(\Zii,\LLmm) = 0$ for $i>0$,
and the Demazure character formula also holds:
$$
\Char \MB_{\ii,\mm} = 
\Lam_{i_1} \om_{i_1}^{m_1} \cdots
\Lam_{i_l} \om_{i_l}^{m_l},
$$
$\om_i$ being the (multplicative) fundamental weights
and $\Lam_i$ the Demazure operators
on the ring of characters of $T$.
\end{prop}
Once again the proof goes through as before,
making appeal to the arguments of \cite{MNW}, Prop 28.

In the case of $G = GL(n)$,
the $\Zii$ for
non-reduced $\ii$ 
again have an explicit interpretation
as configuration varieties.
This is clear from the fiber-product realization
$\Zii \cong \Fii$:
each extra factor in the Bott-Samelson variety
corresponds to one new space in the data of the
configuration variety.  

For example, for $G=GL(3)$ and
$\ii = 2112$, the Bott-Samelson variety is:
$$
\Zii \cong \left\{
\begin{array}{c}
(V_2,V_1,V_1',V_2') \in 
\Gr(2)\times \Gr(1)\times \Gr(1)\times \Gr(2) \\[.1em]
\mbox{ with }\ \CC^1\! \subset\! V_2\! \supset V_1
\ \mbox{  and  }\ 
 V_2\! \supset\! V_1'\! \subset V_2' 
\end{array}
\right\}.
$$

\end{document}